\documentclass[preprint]{elsarticle}
\usepackage{mathrsfs}
\usepackage{amsmath,amssymb}
\usepackage{lineno,hyperref}
\modulolinenumbers[5]

\newtheorem{thm}{Theorem}
 \newtheorem{lem}[thm]{Lemma}
  \newtheorem{cor}[thm]{Corollary}
    \newtheorem{conj}[thm]{Conjecture}
 \newtheorem{defn}[thm]{Definition}
 \newtheorem{prop}[thm]{Proposition}
 \newdefinition{rmk}{Remark}
 \newproof{pf}{Proof}
 \newproof{poti}{Proof of Theorem \ref{thm1}}
 \newproof{potii}{Proof of Theorem \ref{thm2}}

\begin{document}

\begin{frontmatter}

\title{The coset factorization of finite cyclic group}

\author[mymainaddress]{Kevin Zhao}
\ead{zhkw-hebei@163.com}

\address[mymainaddress]{Department of Mathematics, South China normal university, Guangzhou 510631, China}

\begin{abstract}
Let $G$ be a finite cyclic group, written additively, and let $A,\ B$ be nonempty subsets of $G$.
We will say that $G= A+B$ is a \textit{factorization} if for each $g$ in $G$ there are unique elements $a,\ b$ of $G$ such that $g=a+b, \ a\in A, b\in B$.
In particular, if $A$ is a complete set of residues $modulo$ $|A|$,
then we call the factorization a \textit{coset factorization} of $G$.

In this paper, we mainly study a factorization $G= A+B$, where $G$ is a finite cyclic group and $A=[0,n-k-1]\cup\{i_0,i_1,\ldots i_{k-1}\}$
with $|A|=n$ and $n\geq 2k+1$.
We obtain the following conclusion:

If $(i)$ $k\leq 2$ or $(ii)$ The number of distinct prime divisors of $gcd(|A|,|B|)$ is at most $1$ or $(iii)$ $gcd(|A|,|B|)=pq$ with $gcd(pq,\frac{|B|}{gcd(|A|,|B|)})=1$,
then $A$ is a complete set of residues $modulo$ $n$.
\end{abstract}


\end{frontmatter}

\section{Introduction}
The factorizations of finite abelian groups are closely related to some geometric problems.
It was introduced by G. Haj\'{o}s \cite{[GH]} for solving a geometric problem posed by H. Minkowski \cite{[HM]}.

Let $G$ be a finite abelian group, written additively, and let $A_1,\ldots,A_n$ be nonempty subsets of $G$.
If for each $g$ in $G$ there are unique elements $a_1,\ldots,a_n$ of $G$ such that
$$g= a_1 +\ldots+ a_n, \ a_1\in A_1,\ldots,a_n\in A_n,$$
then we say that $G= A_1+\ldots+A_n$ is a \textit{factorization} of $G$ and call $g= a_1 +\ldots+ a_n$ to be the \textit{factorization} of $g$.
If $G=A+B$ is a factorization satisfying that $B$ can be replaceable by a subgroup,
then we call the factorization the \textit{coset factorization} of $G$.
It is easy to see that if $G$ is cyclic, then the definition is equivalent to that
$A$ is a complete set of residues $modulo$ $|A|$.
A subset $A$ of $G$ is called \textit{normalized} if $0\in A$.
A factorization $G= A_1+\ldots+A_n$ is termed normalized if each $A_i$ is normalized.

In this paper, we main investigate the \textit{coset factorization} of a finite cyclic group.
It is closely related to "periodic" and "replaceable".

A nonempty subset $A$ of an abelian group $G$ is defined to be \textit{periodic} if there is an element $g$ of $G$ such that $A + g = A$ and $g\neq 0$.
Let $M=\{g\in G:A + g = A\}$.
It is easy to see that $M$ is a subgroup of $G$.
We call $M$ the \textit{stable subgroup} of $A$ and we can write $A$ as the union of some distinct cosets of $M$
$$A=\cup_{i=1}^{\ell}(b_i+M).$$
Obviously, periodicity of $A$ is a weakening of $A$ being a subgroup of $G$.
A subset $A$ of a group $G$ is said to be \textit{replaceable} by a subset $B$ of $G$ if, whenever $A + C = G$ is a factorization of $G,$ so also is $B + C = G$.
One can find more results about them in \cite{[SS]}.

For a factorization $G = A+B$, one can studying whether it is a coset factorization through the following two aspects:
\begin{description}
  \item[(i)] Study the periodicity of factors to determine whether they are subgroups;
  \item[(ii)] Study whether a factor can be replaceable by a subgroup or a set of coset representative elements of a subgroup.
\end{description}

Furthermore, the factorization of a finite cyclic group is closely related to the splitting problem.
For the relationship one can refer to \cite{[SS], {SS86}, {SS87}}.
Stein \cite{St67} first studied the splitting problem and showed its equivalence to the problem of tiling the Euclidean space by translates of certain polytope composed of unit cubes.

Let $G$ be a finite group, written additively, $M$ a set of integers, and $S$ a subset of $G$.
We will say that $M$ and $S$ form a \textsl{splitting} of $G$ if every nonzero element $g$ of $G$ has a unique representation of the form $g=ms$ with $m\in M$ and $s\in S$, while $0$ has no such representation.
(Here "$ms$" denotes the sum of $m$ $s$'s if $m\geq 0$, and $-((-m)s)$ if $m<0$.)
We will write "$G\setminus \{0\}=MS$" to indicate that $M$ and $S$ form a splitting of $G$.
$M$ will be referred to as the multiplier set and $S$ as the splitting set.
We will also say that $M$ splits $G$ with splitting set $S$, or simply that $M$ splits $G$.

Let $G$ be a cyclic group of order a prime $p$ and let $G\setminus \{0\}=MS$ be a splitting.
Obviously, there exist sets $A$ and $B$ such that
$$G\setminus \{0\}=\mathbb{Z}_p^*=\langle g\rangle=MS=\{g^i:i\in A\}\cdot \{g^j: j\in B\}.$$
That is, $$\mathbb{Z}_{p-1}=A+B$$ is a factorization.
We call $G\setminus \{0\}=MS$ to be a \textit{coset splitting} of $G$, if $S$ can be replaceable by a multiplicative subgroup of $\mathbb{Z}_p^*$ or $M$ is a set of coset representative elements of some multiplicative subgroup.
Thus we say that $M$ coset splits $G$.
In addition, it is easy to see that $\mathbb{Z}_p\setminus \{0\}=MS$ be a coset splitting if and only if $\mathbb{Z}_{p-1}=A+B$ is a coset factorization.
For the coset splitting of $\mathbb{Z}_p$, one can refer to \cite{[GS],S87}.
In \cite{S87}, Stein obtained a sufficient and necessary condition.
In \cite{[GS]}, the authors  called for studying the following problem:

\textsl{Find necessary and sufficient conditions that $M=\{m^{b_1},m^{b_2},\ldots ,m^{b_n}\}$ splits an abelian group.}

They made a preliminary study of this problem and obtained a simple result.

\begin{lem}[\cite{[GS]}, Theorem 5.6] \label{inverse-n-1}
Let $m,$ $n$ be integers with $m\geq 2$ and $n\geq 3$. Suppose $M=\{1,m,m^2,\ldots ,m^{n-2},m^j\}$.
If $M$ splits an abelian group $G$ and $|M|=n$, then $j\equiv n-1$ (mod $n$).
\end{lem}

For this problem, we have the following conjecture:

\begin{conj} \label{mainconjspl}
Let $m,$ $n$, $k$ be positive integers with $m\geq 2$.
Suppose $$M=\{1,m,m^2,\ldots ,m^{n-k-1},m^{i_0},m^{i_1},\ldots ,m^{i_{k-1}}\}.$$
If $M$ splits an abelian group $G$ and $|M|=n\geq 2k+1$, then $\{i_0,i_1,\ldots ,i_{k-1}\}\equiv \{n-k,n-k+1,\ldots ,n-1\}$ (mod $n$).
\end{conj}

Similarly, for the factorizations of finite cyclic groups, we have the following conjecture:
\begin{conj} \label{mainconj}
Let $k$, $\omega$, $n$ be positive integers and let $\mathbb{Z}_{\omega}$ be a finite cyclic group of order $\omega$.
Suppose that $$A=[0,n-k-1]\cup\{i_0,i_1,\ldots i_{k-1}\}$$
is a subset of $\mathbb{Z}_{\omega}$ with $|A|=n\geq 2k+1$.
If $A$ is a direct factor of $\mathbb{Z}_{\omega}$, then $\{i_0,i_1,\ldots ,i_{k-1}\}\equiv \{n-k,n-k+1,\ldots ,n-1\}$ (mod $n$).
\end{conj}

\rmk It is easy to see that if Conjecture \ref{mainconjspl} holds,
then the condition that $n\geq 2k+1$ is optimal.
Let $G=\mathbb{Z}_p$ be a cyclic group of a prime order $p$ and $M=\{1,m,m^2,\ldots ,m^{n-k-1},m^{n},m^{n+1},\ldots ,m^{n+k-1}\}$ with $n=2k\geq 2$.
If $4k|ord_p(m)$,
then we have that
$\langle m\rangle=M\cdot (\langle m^{4k}\rangle\cup m^k\langle m^{4k}\rangle),$
where $\langle m\rangle$ is a subgroup of $\mathbb{Z}_p^*$.
It follows that $M$ splits $\mathbb{Z}_p$ with
the splitting set $\bigcup_{i=1}^ta_i(\langle m^{4k}\rangle\cup m^k\langle m^{4k}\rangle)$,
where $a_1,\ldots ,a_t$ are the coset representations of $\langle m\rangle$ in $\mathbb{Z}_p^*$.
However, $\{n,n+1,\ldots ,n+k-1\}\not\equiv \{n-k,n-k+1,\ldots ,n-1\}$ (mod $n$).
Hence, $n\geq 2k+1$ is optimal for Conjecture \ref{mainconjspl}.

In the following, we show that Conjecture \ref{mainconjspl} is covered by Conjecture \ref{mainconj}.
Thus $n\geq 2k+1$ is also optimal for Conjecture \ref{mainconj}.
For the proof, we will use the following definitions and lemmas.

\begin{defn}[\cite{[H]}, Definition 0.0]\label{def2}
A splitting $G\setminus \{0\}= MS$ of a finite abelian group $G$ is called \textit{nonsingular}
if every element of $M$ is relatively prime to $|G|$; otherwise the splitting is called \textit{singular}.
\end{defn}

\begin{lem}[\cite{[GS]}, Lemma 5.2]\label{relativelyprime}
Let $M$ be a finite set of non-zero integers such that $1\in M$ and $|M|>1$.
Suppose that $q$ is a prime which divides each $m\in M$, $m\neq 1$.
If $M$ splits a group $G$,
then $q$ does not divide $|G|$.
\end{lem}

\begin{lem}[\cite{[H]}, Theorem 2.2.3] \label{nonsingular-sp}
Let $G$ be a finite group and $M$ a set of nonzero integer. Then
$M$ splits $G$ nonsingularly if and only if $M$ splits $\mathbb{Z}_p$ for each prime divisor $p$ of $G$.
\end{lem}

\begin{lem}\label{p-cyclic-group}
Let $\mathbb{Z}_p$ be a cyclic group of a prime order $p$ and $M$ a subset of $\mathbb{Z}_p\setminus \{0\}$.
Then $M$ splits $\mathbb{Z}_p$ if and only if $M$ splits $H$ where $H$ is a subgroup of $\mathbb{Z}_p^*$ with $\langle M\rangle\subseteq H$.

In particular, $M$ splits $\mathbb{Z}_p$ if and only if $M$ splits $\langle M\rangle$.
\end{lem}
\pf
If $M$ splits $\mathbb{Z}_p$ with splitting set $S$, then $H=H\cap \mathbb{Z}_p\setminus \{0\}=H\cap MS$.
Since $s\cdot M\subseteq s\cdot H$ and $s\cdot H\cap H=\phi$ for any $s\in S\setminus H$,
we have $H\subseteq M\cdot (S\cap H)$.
It is easy to see that $M\cdot (S\cap H)\subseteq H\cap MS=H$.
Hence, $H=M\cdot (S\cap H)$ is a splitting.

Suppose that $M$ splits $H$.
Set $H=M\cdot S_1$ and $\mathbb{Z}_p^*=\cup_{i=1}^k a_iH$
where $a_1,\ldots ,a_k$ are the coset representations of $H$ in $\mathbb{Z}_p^*$.
Thus $\mathbb{Z}_p\setminus \{0\}=\mathbb{Z}_p^*=\cup_{i=1}^k a_iH=\cup_{i=1}^k a_i(MS_1)=M\cdot (\cup_{i=1}^k a_iS_1)$. The proof is complete.

\qed

\begin{prop} \label{cover}
If Conjecture \ref{mainconj} holds, then Conjecture \ref{mainconjspl} also holds.
\end{prop}

\pf
In Conjecture \ref{mainconjspl}, Lemma \ref{relativelyprime} implies that $M$ splits $G$ nonsingularly.
By Lemma \ref{nonsingular-sp}, we have that $M$ splits $\mathbb{Z}_p$ for each prime divisor $p$ of $|G|$.
From Lemma \ref{p-cyclic-group} it follows that $M$ splits $\langle M\rangle$.
For $n\geq 2k+1\geq 3$, we have $\langle M\rangle=\langle m\rangle$.
Thus there exists a subset $S$ of $\langle M\rangle$ such that
$$\langle M\rangle=\langle m\rangle=M\cdot S$$ is a splitting.
Hence, $A=[0,n-k-1]\cup\{i_0,i_1,\ldots i_{k-1}\}$ with $|A|=n\geq 2k+1$ is a direct factor of $\mathbb{Z}_{ord_p(m)}$.
By Conjecture \ref{mainconj} we complete the proof.

\qed

Our main results are the following:

\begin{thm}\label{mainthm}
Let $G$ be a cyclic group.
Suppose that $G = A+B$ is a factorization and $A=[0,n-k-1]\cup\{i_0,i_1,\ldots i_{k-1}\}$
with $|A|=n\geq 2k+1$.
If either $\mu(gcd(|A|,|B|))\leq 1$ or $gcd(|A|,|B|)=pq$ with $gcd(pq,\frac{|B|}{gcd(|A|,|B|)})=1$,
where $\mu(gcd(|A|,|B|))$ is the number of distinct prime divisors of $gcd(|A|,|B|)$,
then $\{i_0,i_1,\ldots ,i_{k-1}\}\equiv \{n-k,n-k+1,\ldots ,n-1\}$ (mod $n$).
\end{thm}

\begin{thm} \label{thm5}
Let $k$, $\omega$, $n$ be positive integers and let $\mathbb{Z}_{\omega}$ be a finite cyclic group of order $\omega$.
Suppose that $$A=[0,n-3]\cup\{i,j\}$$
is a subset of $\mathbb{Z}_{\omega}$ with $|A|=n\geq 5$.
If $A$ is a direct factor of $\mathbb{Z}_{\omega}$, then $\{i,j\}\equiv \{n-2,n-1\}$ (mod $n$).
\end{thm}

By Proposition \ref{cover} and Theorem \ref{thm5}, we immediately obtain the following corollary:

\begin{cor} \label{thm5cor}
Let $m,$ $n$ be integers with $m\geq 2$ and $n\geq 5$.
Suppose $M=\{1,m,m^2,\ldots ,m^{n-3},m^i,m^j\}$.
If $M$ splits an abelian group $G$ and $|M|=n$, then $\{i,j\}\equiv \{n-2,n-1\}$ (mod $n$).
\end{cor}

\section{Preliminaries}

In this paper, our notations are coincident with \cite{[GS],[H],SS94,[SS]}
and we briefly present some key concepts.
For real numbers $a,b\in \mathbb{R}$, we set $[a,b] = \{x\in \mathbb{Z} | a\leq x\leq b\}$.
For positive integers $n$ and $g$ with $gcd(n,g)=1$,
let $ord_n(g)$ denote the minimal positive integer $l$ such that $g^l\equiv 1$ (mod $n$).
For a positive integer $n$, let $\mu(n)$ denote the number of distinct prime divisors of $n$,
let $\nu(n)$ denote the number of prime divisors (not necessarily distinct) of $n$.
If $M$ is a subset of a multiplicative group $G$, then $\langle M\rangle$ is a subgroup of $G$ with a generating set $M$.

In order to show our results, we need certain techniques from the character theory.

A character $\chi$ of a finite abelian group $G$ is a mapping from $G$ to the multiplicative group of roots of unity such that $\chi(a + b) = \chi(a)\chi(b)$ for all $a,\ b\in G$.
So $\chi$ is a homomorphism from an additive abelian group to a multiplicative abelian group.
The equations $\chi(0) = \chi(0 + 0) = \chi(0)\chi(0)$ give that $\chi(0) = 1$.
The kernel of $\chi$ is defined by $$Ker\chi = \{a\in G : \chi(a) = 1\}$$ and is a subgroup of $G$.
The unity character maps every element in $G$ to $1$.

If $A$ is a subset of $G$, we define the image under $\chi$ of $A$ by $$\chi(A) =\sum_{a\in A} \chi(a).$$
If $A$ is the empty set, then define $\chi(A)$ to be $0$.
Of course, if $\chi$ is the unity character, then $\chi(A) =|A|$.
We define the annihilator of $A$ by $$Ann(A) = \{\chi : \chi(A) = 0\}.$$
If there is a $g$ in $G$ for which $\chi(g)\neq 1$, then from $\chi(G + g)-\chi(G) = 0$ one obtains that $(\chi(g)-1)\chi(G) = 0$ and thus that $\chi(G) = 0$ for all non-unity characters $\chi$ of $G$.
If $A,\ B$ are subsets of $G$ such that $A + B$ is a direct sum, then $\chi(A + B) = \chi(A)\chi(B)$ for all characters $\chi$ of $G$.

If $g$ is a generator element of $G$ and $|g|= n$, then $$(\chi(g))^n = \chi(ng) = \chi(0) = 1$$ and so $\chi(g)$ is an $(n)$th root of unity.
Conversely, if $\rho$ is an $(n)$th root of unity, then we may define a character $\chi$ of $G$ by setting $\chi(kg) = \rho ^k$.
Thus the number of distinct characters of $G$ is equal to $|G| = n$.

The following lemmas will be used repeatedly:

\begin{lem}[\cite{[SS]}, Theorem 7.1] \label{prime-factor}
Let $m$ and $n$ be relatively prime positive integers.
If $A=\{a_1,\ldots ,a_m\}$ and $B=\{b_1,\ldots ,b_n\}$ are sets of integers such that their sum set
$$A+B=\{a_i+b_j:1\leq i\leq m, 1\leq j\leq n\}$$
is a complete set of representatives $modulo$ $mn$,
then $A$ is a complete set of residues $modulo$ $m$ and $B$ is a complete set of residues $modulo$ $n$.
\end{lem}

\begin{cor} \label{cor-factor-integers}
Let $m,$ $n$ be positive integers and let $p$ be a prime.
Suppose $$M=\{1,m^{i_1},m^{i_2},\ldots ,m^{i_{n-1}}\}.$$
If $M$ splits $\mathbb{Z}_p$, $|M|=n$ and $gcd(\frac{ord_p(m)}{n},n)=1$,
then $\{0,i_1,\ldots ,i_{n-1}\}$ is a complete set of residues modulo $n$.
\end{cor}

\pf
Since $M$ splits $\mathbb{Z}_p$ and $\langle M\rangle\subseteq\langle m\rangle$, by Lemma \ref{p-cyclic-group},
we have that $M$ splits $\langle m\rangle$.
Set $\langle m\rangle=M\cdot S$ with $S=\{1,m^{s_1},\ldots,m^{s_{|S|-1}}\}$.
Thus $$Z_{ord_p(m)}=\{0,i_1,\ldots ,i_{n-1}\}+\{0,s_1,\ldots ,s_{|S|-1}\}$$
is a factorization.
For $(\frac{ord_p(m)}{n},n)=1$, we have $(|M|,|S|)=1$.
By Lemma \ref{prime-factor}, we obtain that $\{0,i_1,\ldots ,i_{n-1}\}$ is a complete set of residues modulo $n$.
\qed

\begin{lem}[\cite{[SS]}, Theorem 7.2] \label{prime-factor-2}
Let $m$ and $n$ be relatively prime positive integers and let $p$ be a prime not dividing $mn$.
If $A$ is a set of $pm$ integers and $B$ is a set of $pn$ integers such that $A+B$ is a complete set of residues modulo $p^2mn$, then either $A$ is a complete set of residues modulo $pm$ or $B$ is a complete set of residues modulo $pn$.
\end{lem}

\begin{lem}[\cite{[SS]}, Theorem 4.2] \label{pq-periodic}
Let $G$ be a cyclic group of order $p^eq^f$, where $p$ and $q$ are distinct primes.
Let $B$ be a subset of order $pq$ such that $\chi(B) = 0$, where $\chi$ is a character of $G$ with $Ker\chi ={0}$.
Then $B$ is periodic.
\end{lem}

\begin{lem}[\cite{[SS]}, Theorem 4.3] \label{periodic}
Let $G$ be a finite cyclic group.
If in a factorization of a group $G$ each factor has either prime power order or order that is the product of two distinct primes, then one factor is periodic.
\end{lem}

\begin{lem}[\cite{[SS]}, Theorem 4.4 and Corollary 4.1] \label{periodic2}
Let $G$ be a finite cyclic group.
If $G = A_1 +\ldots+A_k +B$ is a factorization in which each factor $A_i$ has order a power of a prime $p$,
then one of the factors is periodic.

In particular, in any factorization of a cyclic group of prime power order one factor is always periodic.
\end{lem}

\begin{lem}[\cite{[SS]}, Theorem 4.5] \label{periodic3}
Let $G$ be a cyclic group of order $p^eq$, where $p$ and $q$ are distinct primes.
If $G = A+B$ is a factorization,
then $A$ or $B$ is periodic.
\end{lem}

\begin{lem}[\cite{[SS]}, Theorem 4.6] \label{periodic4}
Let $G$ is a finite cyclic group with $2\leq \nu(|G|)\leq 4$.
If $G = A+B$ is a factorization,
then $A$ or $B$ is periodic.
\end{lem}

\begin{lem}[\cite{[SS]}, Theorem 3.17] \label{k-replace}
If $A$ is a direct factor of a finite abelian group $G$ and $k$ is an integer relatively prime to $|A|$,
then $A$ is replaceable by $kA$.
\end{lem}

\begin{lem}[\cite{[SS]}, Theorem 3.10] \label{fac-equ}
Let $H$, $K$ be subgroups of a finite abelian group $G$ and let $A$ be a non-empty subset of $G$.
The next two statements are equivalent:
\begin{description}
  \item[(i)] $Ann(H)\cap Ann(K)\subseteq Ann(A).$
  \item[(ii)] There exist subsets $B$, $C$ of $G$ such that $A = (H +B)\cup(K +C)$, where the sums $H + B$ and $K + C$ are direct and the union is disjoint.
\end{description}
\end{lem}

\section{The proof of Theorem \ref{mainthm}}

Firstly, a proof of Conjecture \ref{mainconj} can be reduced to a verification of the following case:

\begin{lem} \label{Elemmid}
If Conjecture \ref{mainconj} holds for $\omega=ab$, where $a=|A|$ and $\{p: p|a, \ p \ is \ a \ prime\}=\{p: p|b, \ p \ is \ a \ prime\}$,
then Conjecture \ref{mainconj} holds.
\end{lem}

\pf
Since $A$ is a direct factor of $\mathbb{Z}_{\omega}$, there exists a subset $B$ of $\mathbb{Z}_{\omega}$ such that
$$\mathbb{Z}_{\omega}=A+B$$
is a factorization and normalized.
Assume that $|A|=am$ and $|B|=bn$ satisfying that $gcd(m,n)=1$, $gcd(ab,mn)=1$ and $\{p: p|a, \ p \ is \ a \ prime\}=\{p: p|b, \ p \ is \ a \ prime\}$.
Thus $\omega=abmn$, $gcd(n,|A|)=1$ and $gcd(m,|B|)=1$.
By the condition hypothesis, we can suppose that $mn>1$.
From Lemma \ref{k-replace}, it follows that
$$\mathbb{Z}_{\omega} = A+B=nA + B = A + mB = nA + mB$$ are all factorizations.
It is easy to see that $\mathbb{Z}_{\omega}$ is a direct sum of subgroups $H$, $K$, $L$ of orders $ab$, $m$, $n$, respectively.
Thus we have that $nA\subseteq H +K$ and that $mB\subseteq H +L$.
Let $k\in K$, $\ell\in L$.
Then for $a\in nA$, $b\in mB$ we have that $a + b\in H + k + \ell$ if and only if $a\in H + k, \ b\in H + \ell$.
Hence,
\begin{gather}
(nA\cap (H + k))+(mB\cap (H + \ell))= H + k + \ell   \label{H+k+l-fac}
\end{gather}
is a factorization, that is,
$$((nA-k)\cap H)+((mB-\ell)\cap H)= H$$ is a factorization for any $k\in K$ and $\ell \in L$.
Thus we must have that $|(nA-k)\cap H|$ and $|(mB-\ell)\cap H|$ are both constants.
In addition, let $k=\ell=0$ and we have that
\begin{gather}
(nA\cap H)+(mB\cap H)= H   \label{H-fac}
\end{gather}
is a factorization.
Set $$|(nA-k)\cap H|=x \ and \ |(mB-\ell)\cap H|=y.$$
Since $nA\subseteq H +K$ and $mB\subseteq H +L$,
we have that $nA=\cup_{k\in K}(nA\cap (H+k))$ and $mB=\cup_{\ell\in L}(mB\cap (H+\ell))$.
It follows that $$|nA|=|A|=am=\sum_{k\in K}|nA\cap (H+k)|=\sum_{k\in K}|(nA-k)\cap H|=x|K|=xm$$
and $$|mB|=|B|=bn=\sum_{\ell\in L}|mB\cap (H+\ell)|=\sum_{\ell\in L}|(mB-\ell)\cap H|=y|L|=yn.$$
Hence, $$|(nA-k)\cap H|=x=a \ and \ |(mB-\ell)\cap H|=y=b.$$
Since $H=\langle mn\rangle$ and $A=[0,am-k-1]\cup\{i_0,i_1,\ldots i_{k-1}\}$ with $am\geq 2k+1$,
we have that $$nA\cap H\supseteq \{0, mn, \ldots, \lfloor\frac{am-k-1}{m}\rfloor mn\}.$$
Set $$C=\{0, mn, \ldots, \lfloor\frac{am-k-1}{m}\rfloor mn\}.$$
From $|nA|=am\geq 2k+1$ it follows that
$$|nA\cap H|\geq 2(|nA\cap H|-|C|)+1.$$
Combining the condition hypothesis and (\ref{H-fac}) yields that
$\frac{1}{mn}(nA\cap H)$ is a complete set of residues $modulo$ $a$ and $mB\cap H=\langle amn\rangle=M\cong \mathbb{Z}_b$.
By (\ref{H+k+l-fac}) we have that $$(nA\cap (H + k))+M= H + k$$ is a factorization for any $k\in K$.
Hence, $$nA+M= H + K,$$ since $nA\subseteq H+K$.
It follows that $$nA+M+L= H + K+L=G.$$
Since $M+L$ is a direct sum and $|M+L|=|M|\cdot|L|=bn$,
we have that $M+L\cong \mathbb{Z}_{bn}$, i.e., $M+L=\langle am\rangle$.
Thus $nA$ is a complete set of residues $modulo$ $am$.
For $gcd(n,am)=1$, we have that $A$ is a complete set of residues $modulo$ $am$ and the proof is complete.

\qed

\begin{lem} \label{A-periodic}
Let $G$ be a cyclic group and $A$ is a direct factor of $G$ with $|A|=n$.
If $$A=[0,n-k-1]\cup\{i_0,i_1,\ldots i_{k-1}\}$$ with $2k+1\leq n<|G|$,
then $A$ is not periodic.
\end{lem}

\pf
Let $G$ have order $mn$.
For $n<|G|$, we have $m>1$.
Now suppose that $A$ is periodic.
Thus there exists a nontrivial subgroup $M$ of $G$ such that $$A=\cup_{i=1}^{\ell}(a_i+M),$$
where $|M|=\frac{n}{\ell}>1$.
It follows that $M=\langle m\ell\rangle.$
Since $A=[0,n-k-1]\cup\{i_0,i_1,\ldots i_{k-1}\}$ with $n\geq 2k+1$,
we can claim that $$\ell\geq n-k,$$
which implies that $\ell> \frac{n}{2}$.
This is in contradiction to $\frac{n}{\ell}>1$.
If $m\ell\leq n-k-1$, then $[0,m\ell-1]$ are contained in distinct cosets of $M$, i.e.,
$\ell\geq m\ell$.
It follows that $m=1$, a contradition.
If $m\ell> n-k-1$, then $[0,n-k-1]$ are contained in distinct cosets of $M$,
i.e., $\ell\geq n-k$.
This completes the proof.
\qed

\begin{lem} \label{B-subgroup}
Let $G$ be a finite cyclic group.
Suppose that $G = A+B$ is a factorization and normalized, where $A=[0,n-k-1]\cup\{i_0,i_1,\ldots i_{k-1}\}$
with $|A|=n$ and $n\geq 2k+1$.
If one of the following conditions is satisfied
\begin{description}
  \item[(i)] either $2\leq \nu(|G|)\leq 4$ or $B$ has prime power order,
  \item[(ii)] $B$ is periodic with $|B|=pq$, where $p$, $q$ are distinct primes,
\end{description}
then $A$ is a complete set of residues $modulo$ $n$,
and $B=\langle n\rangle$ is a subgroup.
\end{lem}

\pf
If $|B|=1$, then the lemma is trivial.
Now suppose that $|B|>1$.
Thus Lemma \ref{A-periodic} implies that $A$ is not periodic.
If condition $(i)$ is satisfied, then
combining Lemma \ref{periodic2} and Lemma \ref{periodic4} yields that $B$ is periodic.
Set $$B=\cup_{i=1}^{\ell}(b_i+M),$$
where $M$ is the stable subgroup of $B$ with $|M|>1$.
It follows that $G=A+\{b_1,\ldots,b_{\ell}\}+M$.
Thus $A+\{b_1,\ldots,b_{\ell}\}$ is a complete set of residues $modulo$ $\ell n$, i.e.,
$$A+\{b_1,\ldots,b_{\ell}\}=\mathbb{Z}_{\ell n}$$
is a factorization.
Again, repeat the above reasoning and we have that $\{b_1,\ldots,b_{\ell}\}$ is periodic.
Continue the above discussion and in the end we can obtain that $A$ is a complete set of residues $modulo$ $n$,
and $B$ is a subgroup, i.e., $B=\langle n\rangle$.

Now suppose that condition $(ii)$ is satisfied.
Since $B$ is periodic with $|B|=pq$, by imitating the above proof,
we have that $$A+\{b_1,\ldots,b_{\ell}\}=\mathbb{Z}_{\ell n}$$
is a factorization, where $A$ is not periodic.
It is easy to see that $\ell$ is a prime or $1$.
By Lemma \ref{periodic2}, we have that $\{b_1,\ldots,b_{\ell}\}$ is periodic, which means that it is a subgroup of $\mathbb{Z}_{\ell n}$.
Hence, $A$ is a complete set of residues $modulo$ $n$, and $B=\langle n\rangle$ is a subgroup.

\qed

\begin{lem}\label{Ethm}
Let $G$ be a finite cyclic group.
Suppose that $G = A+B$ is a factorization and $A=[0,n-k-1]\cup\{i_0,i_1,\ldots i_{k-1}\}$
with $|A|=n$ and $n\geq 2k+1$.
If one of the following conditions is satisfied
\begin{description}
  \item[(i)] $G$ has prime power order,
  \item[(ii)] $|G|=p^eq^f$, where $p$, $q$ are distinct primes and $e,\ f$ are positive integers with $e>1$, $f>1$, and $B$ has order that is the product of two distinct primes,
\end{description}
then $\{i_0,i_1,\ldots ,i_{k-1}\}\equiv \{n-k,n-k+1,\ldots ,n-1\}$ (mod $n$).
\end{lem}

\pf
If condition $(i)$ is satisfied, then
Lemma \ref{B-subgroup} completes our proof.
If condition $(ii)$ is satisfied, then
we have that $$|B|=pq\ and \ |A|=p^{e-1}q^{f-1}$$
with $e-1>0$ and $f-1>0$.
In addition, $G$ has subgroups $H,\ K$ of orders $p^e,\ q^f$, respectively, and $G = H + K$ is a direct sum.
Let $H =\langle q^f\rangle= \langle h\rangle,\ K = \langle p^e\rangle= \langle g\rangle$,
$P = \langle p^{e-1}h\rangle= \langle p^{e-1}q^f\rangle$ and $Q = \langle q^{f-1}g\rangle= \langle q^{f-1}p^e\rangle$.
Let $\chi$ be a character of $G$ with $Ker\chi ={0}$.
It is easy to see that $\chi$ is a character of $G$ with $Ker\chi = \{0\}$ if and only if
$\chi(h)$ has order $p^e$ and $\chi(g)$ has order $q^f$ if and only if
$$\chi(p^{e-1}h)\neq 1,\ \chi(q^{f-1}g)\neq 1$$ if and only if
$$\chi(P) = 0,\ \chi(Q) = 0.$$
Since $\chi(G)=\chi(A+B)=\chi(A)\chi(B)=0$,
we have $\chi(A)=0$ or $\chi(B)=0$.
If $\chi(A)=0$, then $\chi'(A)=0$ for any character $\chi'$ of $G$ with $Ker\chi' = \{0\}$.
It follows since $G$ be a cyclic group.
Hence, $$Ann(P)\cap Ann(Q)\subseteq Ann(A).$$
So we may apply Lemma \ref{fac-equ} to conclude that there exist subsets $E$, $F$ such that
\begin{gather}
A = (P + E)\cup(Q + F),   \label{A-bingfac}
\end{gather}
where the sums are direct and the union is disjoint.
It follows that $|A|=p^{e-1}q^{f-1}=|P|\cdot |E|+|Q|\cdot |F|=p|E|+q|F|$, which implies that
$q||E|$, $p||F|$ and $p^{e-2}q^{f-2}=\frac{|E|}{q}+\frac{|F|}{p}$, since $e>1$ and $f>1$.
Set $$|E|=qe_0\ and \ |F|=pf_0.$$
Since $[0,n-k-1]\subseteq A$ with $n=|A|=p^{e-1}q^{f-1}\geq 2k+1$, $P = \langle p^{e-1}q^f\rangle$ and $Q = \langle q^{f-1}p^e\rangle$,
it is easy to see that $n-k-1<p^{e-1}q^f$ and $n-k-1<q^{f-1}p^e$, which implies that $[0,n-k-1]$ are contained in distinct cosets of $P$ and $Q$.
By $(\ref{A-bingfac})$, we have that $$\frac{p^{e-1}q^{f-1}+1}{2}\leq n-k\leq |E|+|F|=qe_0+pf_0=p^{e-2}q^{f-1}+(p-q)f_0.$$
This is impossible.
Hence, we must have that $$\chi(B) = 0.$$
For $|B|=pq$, Lemma \ref{pq-periodic} implies that $B$ is periodic.
By Lemma \ref{B-subgroup}, we complete the proof.

\qed

\textit{The proof of Theorem \ref{mainthm}:}
If $\mu(gcd(|A|,|B|))=0$, then $gcd(|A|,|B|)=1$.
By Lemma \ref{prime-factor}, the theorem holds.
If $\mu(gcd(|A|,|B|))=1$, then combining Lemma \ref{Elemmid} and Lemma \ref{Ethm} $(i)$ yields our conclusion.
If $gcd(|A|,|B|)=pq$ with $gcd(pq,\frac{|B|}{gcd(|A|,|B|)})=1$,
then combining Lemma \ref{Elemmid} and Lemma \ref{Ethm} $(ii)$ yields our conclusion.

\qed

\section{The proofs of Theorem \ref{thm5}}

\textit{The proof of Theorem \ref{thm5}:}
Let $\mathbb{Z}_{\omega}=A+S$ be the factorization and normalized,
where $A=[0,n-3]\cup\{i,j\}$ is a subset of $\mathbb{Z}_{\omega}$ with $|A|=n$ and $n\geq 5$,
and $S=\{0,s_1,\ldots ,s_{|S|-1}\}$.
Thus we can suppose $n-2\leq i<j\leq \omega-1$.
If $i=n-2$, then by Lemma \ref{inverse-n-1}, we get our result.
Now suppose $i>n-2$.

Step $1:$ We will prove $|(j+S)\bigcap [i'-(n-2),i'-1]|\leq 1$ for any $i'\in i+S$.

Suppose $|(j+S)\bigcap [i'-(n-2),i'-1]|\geq 2$.
Then there exist two distinct elements $t_1$ and $t_2$ in $S$ such that $\{j+t_1,\ j+t_2\}\subseteq [i'-(n-2),i'-1]$.
Thus $1\leq |t_1-t_2|\leq n-3$. This is in contradiction to that $\mathbb{Z}_{\omega}=A+S$ is a factorization.

Similarly, we can prove that $|(j+S)\bigcap [i'+1,i'+(n-2)]|\leq 1$, $|(i+S)\bigcap [j'-(n-2),j'-1]|\leq 1$
and $|(i+S)\bigcap [j'+1,j'+(n-2)]|\leq 1$ for any $i'\in i+S$ and $j'\in j+S$.

Step $2:$ We will prove that if $|(j+S)\bigcap [i+1,i+(n-2)]|=1$,
then $(j+S)\bigcap [i+1,i+(n-2)]=\{i+1\}$.

Set $(j+S)\bigcap [i+1,i+(n-2)]=\{i+k\}=\{j+t\}$,
where $k\in [1,n-2]$ and $t\in S$.
Thus $i=j'-k\in [j'-(n-2),j'-1]$ with $j'=j+t$.
Suppose $k>1$.
Thus $j'-(n-2)\leq i<i+1\leq j'-1<i+(n-2).$
It is easy to see that $(i+S)\bigcap[i+1,i+(n-3)]=\phi$ and $(j+S)\bigcap[j'-(n-3),j'-1]=\phi$.
Thus $(i+S\cup j+S)\bigcap [i+1,j'-1]=\phi$, that is $[i+1,j'-1]\subseteq [0,n-3]+S$.
Since $i\in i+S$ and $j'\in j+S$,
there must be some $s\in S$ such that $[0,n-3]+s=[i+1,j'-1]$.
This is a contradiction, since $|[i+1,j'-1]|=j'-1-i<(i+(n-2))-i=n-2=|[0,n-3]+s|$.

Similarly, we can prove that if $|(j+S)\bigcap [i-(n-2),i-1]|=1$,
then $(j+S)\bigcap [i-(n-2),i-1]=\{i-1\}$;
if $|(i+S)\bigcap [j+1,j+(n-2)]|=1$,
then $(i+S)\bigcap [j+1,j+(n-2)]=\{j+1\}$;
if $|(i+S)\bigcap [j-(n-2),j-1]|=1$,
then $(i+S)\bigcap [j-(n-2),j-1]=\{j-1\}$.

Step $3:$In the following, we distinguish several cases:

Case $1:$ $|(j+S)\bigcap [i+1,i+(n-2)]|=1$ and $|(i+S)\bigcap [j-(n-2),j+(n-2)]|=0$.

Thus we can claim that $$\{kn,-i+kn-2,-j+kn-1\}\subseteq S$$ for any positive integer $k$.

By the claim, we have $\langle n\rangle\subseteq S$.
Since $|S|=\frac{\omega}{n}=|\langle n\rangle|$, we have $\langle n\rangle=S$.
Thus $\{i,j\}+S=\{n-2,n-1\}+S$.
It is easy to see that $|(i+S)\bigcap [j-(n-2),j+(n-2)]|=1$, a contradiction.

In the following, we will show that the claim is true.

From Step $2$ it follows that $(j+S)\bigcap [i+1,i+(n-2)]=\{i+1\}$.
Thus there exists $t\in S$ such that $j+t=i+1$, i.e., $i-j+1\in S$.
For $n-3\geq 2$ we have $i+2\not\in i+S\bigcup j+S$ which means $i+2\in [0,n-3]+S$.
Since $i+1\in j+S$, we must have $i+2\in 0+S=S$.
For $i-j+1\in S$ one can easily obtain $i-1\not\in i+S\bigcup j+S$, i.e., $i-1\in [0,n-3]+S$.
It follows that $i-1\in (n-3)+S$, since otherwise $i\in [0,n-3]+S$,
which is in contradiction to $i\in i+S$.
Hence, $i-n+2\in S$.
Since $|(i+S)\bigcap [j-(n-2),j+(n-2)]|=0$,
we have that $j+1\not\in i+S\bigcup j+S$ and $j-1\not\in i+S\bigcup j+S$.
By imitating the above proof, we have $\{j-n+2,j+1\}\subseteq S$.
Therefore, \begin{gather}
\{0,i+2,j+1,i-j+1\}\subseteq S.   \label{shuyuS}
\end{gather}

Since $[0,n-3]+S$ is a direct sum, we have that $[1,n-3]\cap S=\phi$ and $([0,n-3]+s)\cap S=\phi$ for all $s\in \{i+2,j+1,i-j+1\}$.
Hence,
\begin{equation}
\label{jiaoSkong}
   \begin{aligned}
   ([1,n-3]\cup [i+2,i+n-1]\cup [j+1,j+n-2]& \\
   \cup [i-j+1,i-j+n-2])\cap S=\phi.&  \\
   \end{aligned}
  \end{equation}

From the factorization of $j+n-1$ it follows that $|(j+n-1-A)\cap S|=1$, i.e.,
$|([j+2,j+n-2]\cup \{j+n-1,j-i+n-1,n-1\})\cap S|=1$.
By (\ref{jiaoSkong}) we have $[j+2,j+n-2])\cap S=\phi$.
Hence, $$|\{j+n-1,j-i+n-1,n-1\}\cap S|=1.$$
Since $(n-3)+(i+2)=i+(n-1)$ and $i+2\in S$, we have $n-1\not\in S$.
Since $j+1\in S$ and $i+(j+1)=0+(i+j+1)$, we have $i+j+1\not\in S$.
Suppose $j+n-1\in S$.
By the factorization of $i+j+n-1$, we have $|(i+j+n-1-A)\cap S|=1$.
Since $i+j+n-1=i+(j+n-1)$,
we must have that $(i+j+n-1-[0,n-3])\cap S=[i+j+2,i+j+n-1]\cap S=\phi$.
Hence, $[i+j+1,i+j+n-2]\cap S=\phi$.
From (\ref{jiaoSkong}) it follows that $\{i+n-2,j+n-2\}\cap S=\phi$.
By the factorization of $i+j+n-2$,
we have that $|([i+j+1,i+j+n-2]\cup \{i+n-2,j+n-2\})\bigcap S|=1$, a contradiction.
Hence, $j+n-1\not\in S$ and $$j-i+n-1\in S.$$

From the factorization of $j+n$ it follows that $|(j+n-A)\cap S|=1$, i.e.,
$|([j+3,j+n-1]\cup \{j+n,j-i+n,n\})\cap S|=1$.
By (\ref{jiaoSkong}) we have $[j+3,j+n-1]\cap S=\phi$.
Hence, $$|\{j+n,j-i+n,n\}\cap S|=1.$$
Since $j-i+n-1\in S$ and $1+(j-i+n-1)=0+(j-i+n)$, we have $j-i+n\not\in S$.
Since $i+(j+1)=0+(i+j+1),$ $j+(i+2)=0+(i+j+2)$ and $\{j+1,\ i+2\}\subseteq S$, we have $\{i+j+1,\ i+j+2\}\cap S=\phi$.
Suppose $j+n\in S$.
By the factorization of $i+j+n$, we have $|(i+j+n-A)\cap S|=1$.
Since $i+j+n=i+(j+n)$,
we must have $(i+j+n-[0,n-3])\cap S=[i+j+3,i+j+n]\cap S=\phi$.
Hence, $[i+j+1,i+j+n-2]\cap S=\phi$.
From (\ref{jiaoSkong}) it follows that $\{i+n-2,j+n-2\}\cap S=\phi$.
By the factorization of $i+j+n-2$,
we have $|([i+j+1,i+j+n-2]\cup \{i+n-2,j+n-2\})\cap S|=1$, a contradiction.
Hence, $j+n\not\in S$ and $$n\in S.$$

From the factorization of $n-2$ it follows that $|(n-2-A)\cap S|=1$, i.e.,
$|([1,n-3]\cup \{n-2,-i+n-2,-j+n-2\})\cap S|=1$.
By (\ref{jiaoSkong}) we have $[1,n-3]\cap S=\phi$.
Since $n\in S$ and $2+(n-2)=0+n$, we have $n-2\not\in S$.
Since $i-j+1\in S$ and $(n-3)+(i-j+1)=i+(-j+n-2)$, we have $-j+n-2\not\in S$.
Hence, $$-i+n-2\in S.$$

From the factorization of $n-1$ it follows that $|(n-1-A)\cap S|=1$, i.e.,
$|([2,n-2]\cup \{n-1,-i+n-1,-j+n-1\})\cap S|=1$.
By (\ref{jiaoSkong}) we have $[1,n-3]\cap S=\phi$.
For $n-2\not\in S$ we have $[2,n-2]\cap S=\phi$.
Since $n\in S$ and $1+(n-1)=0+n$, we have $n-1\not\in S$.
Since $-i+n-2\in S$ and $1+(-i+n-2)=0+(-i+n-1)$, we have $-i+n-1\not\in S$.
Hence, $$-j+n-1\in S.$$
Thus our claim is true for $k=1$.
We proceed by induction.

Suppose $k>1$ and the claim is true for $k-1$.
Thus $$\{(k-1)n,-i+(k-1)n-2,-j+(k-1)n-1\}\subseteq S.$$
Since $[0,n-3]+S$ is a direct sum, we have that $([0,n-3]+s)\cap (x+S)=\phi$ for all $x\in \{0,i,j\}$ and $s\in \{(k-1)n,-i+(k-1)n-2,-j+(k-1)n-1\}$.
It is easy to see that $(k-1)n+[0,n-3]=[(k-1)n,kn-3]=i+[-i+(k-1)n,-i+kn-3]=j+[-j+(k-1)n,-j+kn-3]$, $-i+(k-1)n-2+[0,n-3]=[-i+(k-1)n-2,-i+kn-5]$ and $-j+(k-1)n-1+[0,n-3]=[-j+(k-1)n-1,-j+kn-4]$.
Hence,
\begin{equation}
\label{jiaoSkongk-1}
   \begin{aligned}
   ([(k-1)n+1,kn-3]\cup [-i+(k-1)n-2,-i+kn-3]& \\
   \cup [-j+(k-1)n-1,-j+kn-3])\cap S=\phi.&  \\
   \end{aligned}
  \end{equation}

From the factorization of $-i+kn-4$ it follows that $|(-i+kn-4-A)\cap S|=1$, i.e.,
$|([-i+(k-1)n-1,-i+kn-4]\cup \{-i-j+kn-4,-2i+kn-4\})\bigcap S|=1$.
By (\ref{jiaoSkongk-1}) we have $[-i+(k-1)n-1,-i+kn-4]\cap S=\phi$.
Since $-j+(k-1)n-1\in S$ and $(n-3)+(-j+(k-1)n-1)=i+(-i-j+kn-4)$, we have $-i-j+kn-4\not\in S$.
Hence, $$-2i+kn-4\in S.$$

Similarly, from the factorization of $-i+kn-3$ it follows that $|([-i+(k-1)n,-i+kn-3]\cup \{-i-j+kn-3,-2i+kn-3\})\bigcap S|=1$.
By (\ref{jiaoSkongk-1}) we have $([-i+(k-1)n-2,-i+kn-3])\cap S=\phi$.
Since $-2i+kn-4\in S$ and $1+(-2i+kn-4)=0+(-2i+kn-3)$, we have $-2i+kn-3\not\in S$.
Hence, $$-i-j+kn-3\in S.$$

By the factorization of $-i+kn-2$,
we have that $|([-i+(k-1)n+1,-i+kn-2]\cup \{-i-j+kn-2,-2i+kn-2\})\bigcap S|=1$.
By (\ref{jiaoSkongk-1}) we have $([-i+(k-1)n-2,-i+kn-3])\cap S=\phi$.
Since $2+(-2i+kn-4)=0+(-2i+kn-2)$, $1+(-i-j+kn-3)=0+(-i-j+kn-2)$ and $\{-2i+kn-4,\ -i-j+kn-3\}\subseteq S$,
we have $\{-i-j+kn-2,-2i+kn-2\}\cap S=\phi$.
Hence, $$-i+kn-2\in S.$$

By the factorization of $-j+kn-2$,
we have $|([-j+(k-1)n+1,-j+kn-2]\cup \{-i-j+kn-2,-2j+kn-2\})\bigcap S|=1$.
Since $i+(-i+kn-2)=j+(-j+kn-2)$, $1+(-i-j+kn-3)=0+(-i-j+kn-2)$ and $\{-i+kn-2,\ -i-j+kn-3\}\subseteq S$,
we have $\{-i-j+kn-2,-j+kn-2\}\cap S=\phi$.
Since $([-j+(k-1)n-1,-j+kn-3])\cap S=\phi$,
we have $$-2j+kn-2\in S.$$

By the factorization of $-j+kn-1$,
we have that $|([-j+(k-1)n+2,-j+kn-1]\cup \{-i-j+kn-1,-2j+kn-1\})\bigcap S|=1$.
Since $1+(-2j+kn-2)=0+(-2j+kn-1)$, $2+(-i-j+kn-3)=0+(-i-j+kn-1)$ and $\{-2j+kn-2,\ -i-j+kn-3\}\subseteq S$,
we have $\{-i-j+kn-1,-2j+kn-1\}\cap S=\phi$.
Since $([-j+(k-1)n-1,-j+kn-3])\cap S=\phi$ and $-j+kn-2\not\in S$,
we have $$-j+kn-1\in S.$$

By the factorization of $kn$,
we have that $|([(k-1)n+3,kn]\cup \{-i+kn,-j+kn\})\bigcap S|=1$.
Since $\{-i+kn-2,$ $-j+kn-1\}\subseteq S$, $i+(-i+kn-2)=0+(kn-2),$ $j+(-j+kn-1)=0+(kn-1)$, $2+(-i+kn-2)=0+(-i+kn)$
and $1+(-j+kn-1)=0+(-j+kn)$,
we have $\{kn-2$, $kn-1$, $-i+kn$, $-j+kn\}\cap S=\phi$.
Since $([(k-1)n+1,kn-3])\cap S=\phi$,
we have $$kn\in S$$ and the proof of the claim is complete.

Case $2:$ $|(j+S)\bigcap [i-(n-2),i-1]|=1$ and $|(i+S)\bigcap [j-(n-2),j+(n-2)]|=0$.

Thus we can claim that $$\{kn,i-j+kn-1,-i+kn-1,-j+kn-2,i+kn+1\}\subseteq S$$ for any positive integer $k$.

By the claim, we have $\langle n\rangle\subseteq S$.
Since $|S|=\frac{\omega}{n}=|\langle n\rangle|$, we have $\langle n\rangle=S$.
Thus $\{i,j\}+S=\{n-2,n-1\}+S$.
It is easy to see that $|(i+S)\bigcap [j-(n-2),j+(n-2)]|=1$, a contradiction.

In the following, we will show that the claim is true.

From Step $2$ it follows that $(j+S)\bigcap [i-(n-2),i-1]=\{i-1\}$.
Thus there exists $t\in S$ such that $j+t=i-1$, i.e., $i-j-1\in S$.
For $n-3\geq 2$ we have $i-2\not\in i+S\bigcup j+S$ which means $i-2\in [0,n-3]+S$.
Since $i-1\in j+S$, we must have $i-2\in (n-3)+S$, i.e., $i-n+1\in S$.
For $i-j-1\in S$ one can easily obtain $i+1\not\in i+S\bigcup j+S$, i.e., $i+1\in [0,n-3]+S$.
It follows that $i+1\in 0+S=S$, since otherwise $i\in [0,n-3]+S$,
which is in contradiction to $i\in i+S$.
Since $|(i+S)\bigcap [j-(n-2),j+(n-2)]|=0$,
we have that $j+1\not\in i+S\bigcup j+S$ and $j-1\not\in i+S\bigcup j+S$.
By imitating the above proof, we have $\{j-n+2,j+1\}\subseteq S$.
Therefore, \begin{gather}
\{0,i+1,j+1,i-j-1\}\subseteq S.   \label{shuyuScase2}
\end{gather}

Since $[0,n-3]+S$ is a direct sum, we have that $[1,n-3]\cap S=\phi$ and $([0,n-3]+s)\cap S=\phi$ for all $s\in \{i+1,j+1,i-j-1\}$.
Hence,
\begin{equation}
\label{jiaoSkongcase2}
   \begin{aligned}
   ([1,n-3]\cup [i+1,i+n-2]\cup [j+1,j+n-2]& \\
   \cup [i-j-1,i-j+n-4])\cap S=\phi.&  \\
   \end{aligned}
  \end{equation}

From the factorization of $i-j+n-3$ it follows that $|(i-j+n-3-A)\cap S|=1$, i.e.,
$|([i-j,i-j+n-4]\cup \{i-j+n-3,-j+n-3,i-2j+n-3\})\cap S|=1$.
By (\ref{jiaoSkongcase2}) we have $[i-j,i-j+n-4]\cap S=\phi$.
Since $(n-4)+(i+1)=j+(i-j+n-3)$ and $i+1\in S$, we have $i-j+n-3\not\in S$.
Since $0\in S$ and $j+(-j+n-3)=0+(n-3)$, we have $-j+n-3\not\in S$.
Hence, $$i-2j+n-3\in S.$$

Similarly, by the factorization of $i-j+n-2$,
we have $|\{i-j+n-2,-j+n-2,i-2j+n-2\}\cap S|=1$.
Since $i-2j+n-3\in S$ and $1+(i-2j+n-3)=0+(i-2j+n-2)$, we have $i-2j+n-2\not\in S$.
Since $i+1\in S$ and $(n-3)+(i+1)=j+(i-j+n-2)$, we have $i-j+n-2\not\in S$.
Hence, $$-j+n-2\in S.$$

By the factorization of $i-j+n-1$,
we have $|\{i-j+n-1,-j+n-1,i-2j+n-1\}\cap S|=1$.
Since $i-2j+n-3\in S$ and $2+(i-2j+n-3)=0+(i-2j+n-1)$, we have $i-2j+n-1\not\in S$.
Since $-j+n-2\in S$ and $1+(-j+n-2)=0+(-j+n-1)$, we have $-j+n-2\not\in S$.
Hence, $$i-j+n-1\in S.$$

By the factorization of $i+n-1$,
we have $|\{i-j+n-1,i+n-1,n-1\}\cap S|=1$.
Since $i-j+n-1\in S$, we have $|\{i+n-1,n-1\}\cap S|=0$.
Since $-j+n-2\in S$ and $j+(-j+n-2)=0+(n-2)$, we have $n-2\not\in S$.
Hence, $|\{n-2,n-1\}\cap S|=0$.
Combining this with the factorization of $n-1$ and (\ref{jiaoSkongcase2}) yields $|\{-i+n-1,-j+n-1\}\cap S|=1$.
Since $-j+n-2\in S$ and $1+(-j+n-2)=0+(-j+n-1)$, we have $-j+n-1\not\in S$.
Hence, $$-i+n-1\in S.$$

By the factorization of $n$,
we have $|\{n,-i+n,-j+n\}\cap S|=1$.
Since $-j+n-2\in S$ and $2+(-j+n-2)=0+(-j+n)$, we have $-j+n\not\in S$.
Since $-i+n-1\in S$ and $1+(-i+n-1)=0+(-i+n)$, we have $-i+n\not\in S$.
Hence, $$n\in S.$$
It follows that $|\{i+n,n+1\}\cap S|=0$.
Hence, $|\{i+n-1,i+n,n+1\}\cap S|=0$.
Combining this with the factorization of $i+n+1$ and (\ref{jiaoSkongcase2}) yields $|\{i+n+1,i-j+n+1\}\cap S|=1$.
Since $i-j+n-1\in S$ and $2+(i-j+n-1)=0+(i-j+n+1)$, we have $i-j+n+1\not\in S$.
Hence, $$i+n+1\in S.$$
Thus our claim is true for $k=1$.
We proceed by induction.

Suppose $k>1$ and the claim is true for $k-1$.
Thus $$\{(k-1)n,i-j+(k-1)n-1,-i+(k-1)n-1,-j+(k-1)n-2,i+(k-1)n+1\}\subseteq S.$$

Since $[0,n-3]+S$ is a direct sum, we have that $([0,n-3]+s)\cap (x+S)=\phi$ for all $x\in \{0,i,j\}$ and $s\in \{(k-1)n,i-j+(k-1)n-1,-i+(k-1)n-1,-j+(k-1)n-2,i+(k-1)n+1\}$.
It is easy to see that $(k-1)n+[0,n-3]=[(k-1)n,kn-3]=i+[-i+(k-1)n,-i+kn-3]=j+[-j+(k-1)n,-j+kn-3]$, $i-j+(k-1)n-1+[0,n-3]=[i-j+(k-1)n-1,i-j+kn-4]=i+[-j+(k-1)n-1,-j+kn-4]$,
$-i+(k-1)n-1+[0,n-3]=[-i+(k-1)n-1,-i+kn-4]$,
$-j+(k-1)n-2+[0,n-3]=[-j+(k-1)n-2,-j+kn-5]$ and
$i+(k-1)n+1+[0,n-3]=[i+(k-1)n+1,i+kn-2]=i+[(k-1)n+1,kn-2]$.
Hence,
\begin{equation}
\label{jiaoSkongk-1case2}
   \begin{aligned}
   &([(k-1)n+1,kn-2]\cup [i-j+(k-1)n-1,i-j+kn-4] \\
   &\cup [-i+(k-1)n-1,-i+kn-3]\cup [-j+(k-1)n-1,-j+kn-3]  \\
   &\cup [i+(k-1)n+1,i+kn-2])\cap S=\phi.  \\
   \end{aligned}
  \end{equation}

By the factorization $i-j+kn-3$,
we have $|\{i-j+kn-3,-j+kn-3,i-2j+kn-3\}\cap S|=1$.
Since $(n-4)+(i+(k-1)n+1)=j+(i-j+kn-3)$ and $i+(k-1)n+1\in S$, we have $i-j+kn-3\not\in S$.
Since $n-3\in S$ and $j+(-j+kn-3)=(n-3)+(k-1)n$, we have $-j+kn-3\not\in S$.
Hence, $$i-2j+kn-3\in S.$$

By the factorization of $i-j+kn-2$,
we have $|\{i-j+kn-2,-j+kn-2,i-2j+kn-2\}\cap S|=1$.
Since $i-2j+kn-3\in S$ and $1+(i-2j+kn-3)=0+(i-2j+kn-2)$, we have $i-2j+kn-2\not\in S$.
Since $i+1\in S$ and $(n-3)+(i+(k-1)n+1)=j+(i-j+kn-2)$, we have $i-j+kn-2\not\in S$.
Hence, $$-j+kn-2\in S.$$

By the factorization of $i-j+kn-1$,
we have $|\{i-j+kn-1,-j+kn-1,i-2j+kn-1\}\cap S|=1$.
Since $i-2j+kn-3\in S$ and $2+(i-2j+kn-3)=0+(i-2j+kn-1)$, we have $i-2j+kn-1\not\in S$.
Since $-j+kn-2\in S$ and $1+(-j+kn-2)=0+(-j+kn-1)$, we have $-j+kn-2\not\in S$.
Hence, $$i-j+kn-1\in S.$$

By the factorization of $i+kn-1$,
we have $|\{i-j+kn-1,i+kn-1,kn-1\}\cap S|=1$.
Since $i-j+kn-1\in S$, we have $|\{i+kn-1,kn-1\}\cap S|=0$.
Since $-j+kn-2\in S$ and $j+(-j+kn-2)=0+(kn-2)$, we have $kn-2\not\in S$.
Hence, $|\{n-2,n-1\}\cap S|=0$.
Combining this with the factorization of $kn-1$ and (\ref{jiaoSkongk-1case2}) yields $|\{-i+kn-1,-j+kn-1\}\cap S|=1$.
Since $-j+kn-2\in S$ and $1+(-j+kn-2)=0+(-j+kn-1)$, we have $-j+kn-1\not\in S$.
Hence, $$-i+kn-1\in S.$$

By the factorization of $kn$,
we have $|\{kn,-i+kn,-j+kn\}\cap S|=1$.
Since $-j+kn-2\in S$ and $2+(-j+kn-2)=0+(-j+kn)$, we have $-j+kn\not\in S$.
Since $-i+kn-1\in S$ and $1+(-i+kn-1)=0+(-i+kn)$, we have $-i+kn\not\in S$.
Hence, $$kn\in S.$$
It follows that $|\{i+kn,kn+1\}\cap S|=0$.
Hence, $|\{i+kn-1,i+kn,kn+1\}\cap S|=0$.
Combining this with the factorization of $i+kn+1$ and (\ref{jiaoSkongk-1case2}) yields $|\{i+kn+1,i-j+kn+1\}\cap S|=1$.
Since $i-j+kn-1\in S$ and $2+(i-j+kn-1)=0+(i-j+kn+1)$, we have $i-j+kn+1\not\in S$.
Hence, $$i+kn+1\in S$$
and the proof of the claim is complete.

Case $3:$ $|(i+S)\bigcap [j-(n-2),j-1]|=1$ and $|(j+S)\bigcap [i-(n-2),i+(n-2)]|=0$.

It suffices to exchange $i$ and $j$ in Case $2$.

Case $4:$ $|(i+S)\bigcap [j+1,j+(n-2)]|=1$ and $|(j+S)\bigcap [i-(n-2),i+(n-2)]|=0$.

It suffices to exchange $i$ and $j$ in Case $1$.

Case $5:$ Either $|(j+S)\bigcap [i+1,i+(n-2)]|=1$, $|(i+S)\bigcap [j-(n-2),j-1]|=1$
or $|(j+S)\bigcap [i-(n-2),i-1]|=1$, $|(i+S)\bigcap [j+1,j+(n-2)]|=1$.

By imitating the proof of Case $1$, we can show that our theorem is true.
Remark that (\ref{shuyuS}) will be replaced by $\{0,i+2,j+1,i-j+1\}\subseteq S$ in the two cases.

Case $6:$ $|(j+S)\bigcap [i-(n-2),i-1]|=1$ and $|(i+S)\bigcap [j-(n-2),j-1]|=1$.

By imitating the proof of Case $2$, we can show that our theorem is true.
Remark that (\ref{shuyuScase2}) will be replaced by $\{0,i+1,j+1,i-j-1\}\subseteq S$ in this case.

Case $7:$ $|(j+S)\bigcap [i+1,i+(n-2)]|=1$ and $|(i+S)\bigcap [j+1,j+(n-2)]|=1$.

By imitating the proof of (\ref{shuyuScase2}) in Case $2$, we have $\{0,i+2,j+2,i-j+1,j-i+1\}\subseteq S$.
From the factorization of $n-2$ it follows that $|(n-2-A)\cap S|=1$, i.e.,
$|([1,n-3]\cup \{n-2,-i+n-2,-j+n-2\})\cap S|=1$.
It is easy to see that $[1,n-3]\cap S=\phi$, $i+(n-2)=(n-4)+(i+2)$, $j+(-i+n-2)=(n-3)+(j-i+1)$ and $i+(-j+n-2)=(n-3)+(i-j+1)$.
This is a contradiction.

Case $8:$ $|(i+S)\bigcap [j-(n-2),j+(n-2)]|=0$ and $|(j+S)\bigcap [i-(n-2),i+(n-2)]|=0$.

Since $|(i+S)\bigcap [j-(n-2),j+(n-2)]|=0$,
we have that $j+1\not\in i+S\bigcup j+S$.
From the proof of Case $1$, it is easy to see that $j+1\in S$.
Repeat the reasoning and from $|(j+S)\bigcap [i-(n-2),i+(n-2)]|=0$ it follows that $i+1\in S$.
Thus $i+j+1\in (i+S)\bigcap(j+S)$, a contradiction.

\qed

\section{Open Problems}

For Lemma \ref{inverse-n-1}, we have the following generation:

\begin{lem} \label{inverse-n-1-genfac}
Let $\omega,$ $n$ be positive integers, $n\geq 3$ and let $\mathbb{Z}_{\omega}$ be a cyclic group of order $\omega$.
Suppose that $A=[0,n-1]\setminus \{i\}\cup \{j\}$ is a subset of $\mathbb{Z}_{\omega}$ with $|A|=n$ and $i\in [0,n-1]$.
If $A$ is a direct factor of $\mathbb{Z}_{\omega}$, then $j\equiv i$ (mod $n$).
\end{lem}

\pf
Let $\mathbb{Z}_{\omega}=A+B$ be the factorization and normalized.
It suffices to show that $B=\langle n\rangle$, i.e., $kn\in B$ for all nonnegative integers $k$.
The proof is by induction on $k$.
It is easy to see that the result is true for $k=0$.
Suppose that $k\geq 1$ and we have proved the result for $k-1$.
By the induction hypothesis we have $$(k-1)n\in B.$$

If $i=n-1$, then $A=[0,n-2]\cup \{j\}$ with $n\geq 3$.
By the factorization of $kn-1$, we have that $$|(kn-1-A)\cap B|=|([(k-1)n+1,kn-1]\cup \{kn-1-j\})\cap B|=1.$$
Since $(k-1)n\in B$ and $0+((k-1)n+t)=t+(k-1)n$ for any $t\in [1,n-2]$, we have that $[(k-1)n+1,kn-2]\cap B=\phi$, i.e, $$exact \ one \ of\ \{kn-1,\ kn-1-j\} \ is \ contained \ in \ B.$$
If $kn-1\in B$, then $kn-1-j\not\in B$.
From the factorization of $kn-1-j$ it follows that $$|(kn-1-j-A)\cap B|=|([(k-1)n+1-j,kn-1-j]\cup \{kn-1-2j\})\cap B|=1.$$
Since $kn-1-j\not\in B$ and $j+((k-1)n+t-j)=t+(k-1)n$ for any $t\in [1,n-2]$, we have $[(k-1)n+1-j,kn-1-j]\cap B=\phi$, i.e., $kn-1-2j\in B$.
Since $kn=1+(kn-1)=j+(kn-j)$ with $kn-1\in B$ and $\{1,\ j\}\subseteq A$,
we have $kn-j\not\in B$.
From the factorization of $kn-j$ it follows that $$|(kn-j-A)\cap B|=|([(k-1)n+2-j,kn-j]\cup \{kn-2j\})\cap B|=1.$$
Since $[(k-1)n+2-j,kn-1-j]\cap B=\phi$ and $kn-j\not\in B$, we have $kn-2j\in B$.
This is a contradiction, since $1+(kn-1-2j)=0+(kn-2j)$ with $\{kn-1-2j,\ kn-2j\}\subseteq B$.
Hence, $$kn-1\not\in B\ and \ kn-1-j\in B.$$
Since $0+(kn-j)=1+(kn-1-j)$ with $kn-1-j\in B$, we have $kn-j\not\in B$.
By the factorization of $kn$, we have $$|(kn-A)\cap B|=|([(k-1)n+2,kn]\cup \{kn-j\})\cap B|=1.$$
Combining the above yields that $$kn\in B$$ and we complete the proof for $i=n-1$.

If $i=0$, then $\mathbb{Z}_{\omega}=(A-1)+B$ is a factorization where $A-1=[0,n-2]\cup \{j-1\}$.
By the proof of $i=n-1$, we obtain our result.

Now suppose $i\in [1,n-2]$.
From the factorization of $(k-1)n+i$ it follows that $$|((k-1)n+i-A)\cap B|=|([(k-2)n+i+1,(k-1)n+i]\setminus \{(k-1)n\}\cup \{(k-1)n+i-j\})\cap B|=1.$$
That is $$|(\{(k-1)n+t:t\in [i-n+1,i]^*\}\cup \{(k-1)n+i-j\})\cap B|=1$$
with $[i-n+1,i]^*=[i-n+1,i]\setminus \{0\}\subseteq [-(n-2),n-2]$.
Since $(k-1)n\in B$, $(-t)+((k-1)n+t)=0+(k-1)n$ for any $t\in [i-n+1,-1]\setminus \{-i\}$,
$0+((k-1)n+t)=t+(k-1)n$ for any $t\in [1,i-1]$,
$(i+1)+((k-1)n-i)=1+(k-1)n$
and $1+((k-1)n+i)=(i+1)+(k-1)n$ with $\{1,i+1\}\subseteq A$,
we have $([(k-2)n+i+1,(k-1)n+i]\setminus \{(k-1)n\})\cap B=\phi$, i.e.,
$$(k-1)n+i-j\in B.$$
From the factorization of $kn$ it follows that $$|(kn-A)\cap B|=|([(k-1)n+1,kn]\setminus \{kn-i\}\cup \{kn-j\})\cap B|=1.$$
Since $(k-1)n\in B$, $0+((k-1)n+t)=t+(k-1)n$ for any $t\in [1,n-1]\setminus \{n-i,i\}$
and $1+((k-1)n+i)=(i+1)+(k-1)n$ with $\{1,i+1\}\subseteq A$,
we have $([(k-1)n+1,kn-1]\setminus \{kn-i\})\cap B=\phi$, that is,
$$exact \ one \ of\ \{kn, \ kn-j\} \ is \ contained \ in \ B.$$
If $n-i\neq i$, then we have $0+(kn-j)=(n-i)+((k-1)n+i-j)$ with $n-i\in [2,n-1]\setminus \{i\}\subseteq A$ and $(k-1)n+i-j\in B$.
If $n-i=i$, then we have $1+(kn-j)=(n-i+1)+((k-1)n+i-j)$ with $\{1,n-i+1\}=\{1,i+1\}\subseteq A$ and $(k-1)n+i-j\in B$.
Hence, $kn-j\not\in B$.
That is $$kn\in B$$ and we complete the proof.

\qed

Repeat the reasoning of Proposition \ref{cover} and one can easily obtain the following result:

\begin{cor} \label{inverse-n-1-gen}
Let $m,$ $n$ be integers with $m\geq 2$ and $n\geq 3$.
Suppose that $M=\{1,m,m^2,\ldots ,m^{n-2},m^{n-1}\}\setminus \{m^i\}\cup \{m^j\}$ with $i\in [0,n-1]$.
If $M$ splits an abelian group $G$ and $|M|=n$, then $j\equiv i$ (mod $n$).
\end{cor}

Hence, we have the following problems:

\begin{description}
  \item[(1)] Let $k$, $\omega$, $n$ be positive integers and let $\mathbb{Z}_{\omega}$ be a finite cyclic group of order $\omega$.
If $A$ is a direct factor of $\mathbb{Z}_{\omega}$ with $|A\cap[0,n-1]|\geq\frac{n+1}{2}$ and $|A|=n$, then  investigate whether $A$ is a complete set of residues modulo $n$.
  \item[(2)] Let $k$, $\omega$, $n$ be positive integers and let $\mathbb{Z}_{\omega}$ be a finite cyclic group of order $\omega$.
If $A$ is a direct factor of $\mathbb{Z}_{\omega}$ satisfying that $|A\cap B|\geq\frac{n+1}{2}$ and $|A|=n$,
where $B$ is a complete set of residues modulo $n$, then investigate whether $A$ is a complete set of residues modulo $n$.
  \item[(3)] Are Conjecture \ref{mainconjspl} and Conjecture \ref{mainconj} equivalent?
\end{description}

\rmk
By imitating the proof of Proposition \ref{cover},
one can easily show that in the following conjecture,
$(i)$ is covered by $(ii)$, and $(i)$, $(ii)$ is equivalent to Conjecture \ref{mainconjspl}, \ref{mainconj}, respectively.
In addition, by Lemma \ref{inverse-n-1}, \ref{inverse-n-1-genfac}, Theorem \ref{thm5} and Corollary \ref{thm5cor},
it is easy to see that the following conjecture is true for $k\leq 2$.

\begin{conj} \label{equal-conj}

\begin{description}
  \item[(i)] Let $m,$ $n$, $k$ be positive integers with $m\geq 2$.
Suppose $$M=\{1,m,m^2,\ldots ,m^{n-k-1},m^{i_0},m^{i_1},\ldots ,m^{i_{k-1}}\}$$ with $n\geq 2k+1$.
Then $M$ splits an abelian group $G$ with $|M|=n$ if and only if for any prime $p||G|$, $\{i_0,i_1,\ldots ,i_{k-1}\}\equiv \{n-k,n-k+1,\ldots ,n-1\}$ (mod $n$) and $n|ord_p(m)$.
  \item[(ii)] Let $k$, $\omega$, $n$ be positive integers and let $\mathbb{Z}_{\omega}$ be a finite cyclic group of order $\omega$.
Suppose that $$A=[0,n-k-1]\cup\{i_0,i_1,\ldots i_{k-1}\}$$
is a subset of $\mathbb{Z}_{\omega}$ with $|A|=n\geq 2k+1$.
Then $A$ is a direct factor of $\mathbb{Z}_{\omega}$ if and only $\{i_0,i_1,\ldots ,i_{k-1}\}\equiv \{n-k,n-k+1,\ldots ,n-1\}$ (mod $n$) and $n|\omega$.
\end{description}

\end{conj}

\section*{Acknowledgments}
This work is supported by NSF of China  (Grant No. 11671153).
The author is sincerely grateful to professor Pingzhi Yuan for his guidance and the anonymous referee
for useful comments and suggestions.

\section*{References}


\end{document}